\documentclass[12pt]{amsart}
\usepackage{amsfonts,latexsym,amsthm,amssymb}
\usepackage[all]{xy}

\input epsf

\newtheorem{theorem}{Theorem}[section]
\newtheorem{lemma}[theorem]{Lemma}
\newtheorem{corol}[theorem]{Corollary}
\newtheorem{prop}[theorem]{Proposition}

{\theoremstyle{definition} \newtheorem{defin}[theorem]{Definition}}
{\theoremstyle{remark} 
\newtheorem{example}[theorem]{Example}}

\newcommand{\Abb}{{\mathbb{A}}}
\newcommand{\Cbb}{{\mathbb{C}}}

\newcommand{\Pbb}{{\mathbb{P}}}
\newcommand{\Zbb}{{\mathbb{Z}}}

\DeclareMathOperator{\Spec}{Spec}
\DeclareMathOperator{\Eu}{Eu}

\newcommand{\x}{$\times$}
\newcommand{\one}{1\hskip-3.5pt1}
\newcommand{\csm}{c_{\text{\rm SM}}}
\newcommand{\cma}{c_{\text{\rm Ma}}}
\newcommand{\ssm}{s_{\text{\rm SM}}}
\newcommand{\sma}{s_{\text{\rm Ma}}}
\newcommand{\ccsm}{\check c_{\text{\rm SM}}}
\newcommand{\ccma}{\check c_{\text{\rm Ma}}}
\newcommand{\cssm}{\check s_{\text{\rm SM}}}
\newcommand{\csma}{\check s_{\text{\rm Ma}}}
\newcommand{\Ch}{\text{\rm Ch}}
\newcommand{\Proj}{\text{\rm Proj}}
\newcommand{\Sym}{\text{\rm Sym}}
\newcommand{\qSym}{\text{\rm qSym}}
\newcommand{\Rees}{\text{\rm Rees}}
\newcommand{\Tors}{\text{\rm Tors}}
\newcommand{\Disc}{\text{\rm Disc}}

\newcommand{\cHom}{\mathcal Hom}
\newcommand{\qBl}{qB\ell}

\title{Shadows of blow-up algebras}
\author{Paolo Aluffi}
\address{Max-Planck-Institut f\"ur Mathematik, Bonn, Germany}
\address{Florida State University, Tallahassee, Florida}

\makeindex

\begin{document}

\begin{abstract}
We study different notions of blow-up of a scheme $X$ along a
subscheme~$Y$, depending on the datum of an embedding of $X$ into an
ambient scheme.  The two extremes in this theory are the ordinary
blow-up, $B\ell_YX$, corresponding to the identity $X\overset
=\hookrightarrow X$, and the `quasi-symmetric blow-up', $\qBl_YX$,
corresponding to an embedding $X\hookrightarrow M$ into a nonsingular
variety $M$. We prove that this latter blow-up is intrinsic of $Y$ and
$X$, and is universal with respect to the requirement of being
embedded as a subscheme of the ordinary blow-up of some ambient space
along~$Y$.

We consider these notions in the context of the theory of characteristic 
classes of singular varieties. We prove that if $X$ is a hypersurface in a 
nonsingular variety and $Y$ is its `singularity subscheme', these two 
extremes embody respectively the {\em conormal\/} and {\em characteristic\/} 
cycles of $X$. Consequently, the first carries the essential information 
computing Chern-Mather classes, and the second is likewise a carrier for 
Chern-Schwartz-MacPherson classes. In our approach, these classes are 
obtained from Segre class-like invariants, in precisely the same way as 
other intrinsic characteristic classes such as those proposed by William 
Fulton, \cite{MR85k:14004}, and by W.~Fulton and Kent Johnson, 
\cite{MR82c:14017}.

We also identify a condition on the singularities of a hypersurface under 
which the quasi-symmetric blow-up is simply the linear fiber space 
associated with a coherent sheaf.
\end{abstract}

\maketitle


\section{Introduction}\label{intro}

It is not hard to see that the {\em conormal\/} cycle of a hypersurface 
$X$ of a nonsingular algebraic variety $M$ can be realized as the cycle 
of the blow-up of $X$ along its singularity subscheme (defined by the 
partials of an equation defining $X$). Our guiding question in this paper 
is, what kind of `blow-up' realizes similarly the much subtler {\em 
characteristic\/} cycle of a hypersurface? We answer this question, and 
extract from our construction a unified approach to different 
characteristic classes associated with a possibly singular hypersurface of 
a nonsingular variety.

The ordinary blow-up of a scheme $X$ along a subscheme $Y$---that is, 
the Proj of the Rees algebra of the ideal sheaf $\mathcal J_{Y,X}$ of $Y$ 
in $X$---has the remarkable property that it can be recovered from the blow-up 
of any ambient scheme $M$ along $Y$, by taking the proper transform of $X$. 
As there are other notions of blow-up, obtained by taking the Proj of 
other `blow-up algebras' (such as the symmetric algebra of $\mathcal 
J_{Y,X}$), it is natural to ask whether there is a `largest' blow-up of $X$ 
along $Y$ that can be embedded in {\em some\/} (ordinary) blow-up of an 
ambient scheme $M$ along $Y$.

In the first part of this paper we construct such a blow-up: we define a 
new {\em quasi-symmetric algebra\/} of an ideal $\mathcal J_{Y,X}$, and 
show that it satisfies the universal property summarized above. In fact, we 
define (Definition~\ref{qsymdef}) {\em a\/} quasi-symmetric algebra for every 
embedding $X\subset M$, then show (Theorem~\ref{concom}) that the limit of 
the corresponding inverse system of algebras equals the quasi-symmetric 
algebra arising for any nonsingular $M$ (otherwise independently of $M$). 
We name the corresponding blow-up the {\em quasi-symmetric blow-up\/} of $X$ 
along $Y$, $\qBl_YX$. We also show (Theorem~\ref{printran}) that this new 
blow-up can be obtained by taking a `principal' transform of $X$ in $B\ell_YM$, 
for any nonsingular variety $M$ containing $X$.

The ordinary Rees blow-up and the new quasi-symmetric blow-up are two 
extremes in a range. In the second part of the paper we consider the case 
in which $X$ is a hypersurface in a nonsingular ambient variety $M$, and we 
take $Y$ to be its {\em singularity subscheme.\/} We find that the two 
extremes live naturally in the projectivized cotangent bundle of $M$, and 
their cycles yield concrete realizations of the {\em conormal,\/} 
resp.~{\em characteristic\/} cycles of $X$ (Theorems~\ref{rees} 
and~\ref{qSym}). As mentioned above, the first of these facts is old fare; 
the second appears to be new, at least in the form given here. Every 
quasi-symmetric blow-up in the range should correspond to a lagrangian cycle 
in the projectivized cotangent bundle; that is, every embedding of $X$ in 
another scheme should determine a constructible function on $X$ by this 
construction. One way to summarize the main results in \S\ref{cycles} is by 
saying that our construction associates the identity $X\overset 
=\hookrightarrow X$ with the Euler obstruction of $X$, and any inclusion 
$X\subset M$ into a nonsingular variety with the constant function $\one_X$.

From the point of view of characteristic classes of singular
hypersurfaces, this means that `Rees is to Mather as quasi-symmetric
is to Schwartz-MacPherson'. In the third part of the paper we show
(Theorem~\ref{shadtheorem}) how to obtain these classes rather
directly from the corresponding blow-up algebras, by a standard
intersection-theoretic operation (which is the `shadow' in the title,
Definition~\ref{shaddef}).  This set-up gives a unified approach---for
hypersurfaces---for the theory of Chern-Mather and
Chern-Schwartz-MacPherson classes together with other intrinsic
classes defined for singular varieties---notably the classes defined
by William Fulton and Kent Johnson in \cite{MR82c:14017}, and those
defined by W.~Fulton in \cite{MR85k:14004}, Example~4.2.6.

We also discuss briefly (\S\ref{xfiles}) an intriguing condition on
the singularities of a hypersurface, under which the quasi-symmetric
algebra of the singularity subscheme equals the symmetric algebra; in
other words, in this case the characteristic cycle of $X$ is the
linear fiber space of the coherent sheaf $\mathcal J_{Y,X}$, and the
Chern-Schwartz-MacPherson class of $X$ can be computed from the
ordinary Segre class of a coherent sheaf.  We point out that this
condition is automatically verified in several standard situations,
and mention an interpretation of the condition in terms of extending
vector fields along pieces of a Whitney stratification of the
hypersurface.
\vskip 6pt

One should wonder whether an intrinsic realization of the
characteristic cycle can be given for more general schemes than
hypersurfaces of nonsingular varieties (as we do here). In the end,
our attention is directed to a coherent sheaf that is present
regardless of whether $X$ is a hypersurface: the cokernel of the dual
of the map on differentials determined by the embedding in a
nonsingular variety. If $X$ is a hypersurface then a quasi-symmetric
algebra can be defined for this sheaf, and our main result shows that
this algebra leads to the Chern-Schwartz-MacPherson class of $X$
(Theorem~\ref{lean}).

This suggests what the shape of an analogous result for arbitrary
schemes might be, but the difficulty in establishing such a general
result should not be underestimated. Indeed, the key technical fact
allowing us to obtain the result for hypersurfaces in this paper
amounts to a specific result relating Fulton-Johnson's classes and
Chern-Schwartz-MacPherson classes of hypersurfaces. This relation has
now been known for the better part of a decade, and studied intensely
from many different viewpoints (cf.~\cite{MR96d:14004},
\cite{MR97m:14003}, \cite{MR2001g:32062},
\cite{MR2001k:14014}, \cite{MR2001i:14009}, \cite{MR2001d:14008}, 
\cite{MR2002a:14004}, \cite{MR1819626}, \cite{MR1795550} and the recent 
\cite{schuer2} to name a few), yet a generalization to arbitrary schemes has 
proved exceedingly elusive. A full analog of the results in this paper
to arbitrary schemes would amount to a solution of this problem.

Our motivation in pursuing this program is twofold. First, we believe that 
it would be highly worthwhile to uncover any functoriality feature of
classes such as Fulton's or Fulton-Johnson's. Chern-Schwartz-MacPherson's 
classes owe their existence precisely to their excellent functoriality 
properties; if such functoriality could be transferred to Segre classes (via 
formulas such as the ones presented in this article), this would offer a new
handle on computing Segre classes, arguably one of the most basic
invariants in intersection theory. Second, formulas such as the ones 
obtained in this paper can be implemented into algorithms running in 
symbolic computation programs such as {\tt Macaulay2} (\cite{M2}). The 
only algorithm known to us for such computations
(\cite{math.AG/0204230}) is woefully slow, and we hope that the
approach presented in this paper may lead to substantially improved
algorithms.
\vskip 6pt

\noindent{\bf Acknowledgments} This work was performed while the author
was visiting the Max-Planck-Institut f\"ur Mathematik in Bonn,
Germany; the present version is a thorough reworking of an MPI 
preprint by the same title. Thanks are due to the MPI for support and for 
the congenial atmosphere, and to Prof.~Marcolli for countless insightful
discussions.


\section{quasi-Symmetric algebras and blow-ups}\label{qsymblowups}

\subsection{}\label{definqsym}
Our rings will be Noetherian, commutative, with~1. Homomorphisms of 
algebras endowed of a natural grading are implicitly understood to 
preserve the grading.

Let $A$ be a ring, and $\mathfrak a$ an ideal of $A$. Let $R$ be a ring 
surjecting onto $A$, and denote by $I$ the inverse image of $\mathfrak a$ in 
$R$. Note that the symmetric algebra $\Sym_R(I)$ maps to both the Rees 
algebra $\Rees_R(I)$ and (by functoriality of $\Sym$) to $\Sym_A(\mathfrak 
a)$.

\begin{defin}\label{qsym}
The {\em quasi-symmetric algebra\/} $\qSym_{R\to A}(\mathfrak a)$ is 
defined by
$$\qSym_{R\to A}(\mathfrak a):=\Sym_A(\mathfrak a)\otimes_{\Sym_R(I)} 
\Rees_R(I)\quad.$$
\end{defin}

A particular case of this notion will be the affine version of our main 
blow-up algebra, cf.~Definition~\ref{absqsym} below. Note that the algebra 
corresponding to the identity is the ordinary Rees algebra:
$$\qSym_{A\to A}(\mathfrak a)=\Rees_A(\mathfrak a)\quad;$$
thus, the ordinary blow-up can be recovered in terms of the operation 
studied here. We will be especially interested in the case corresponding 
to epimorphisms $R\to A$ with $R$ suitably `nice'; we begin by recording a 
few properties of the local version of the more general notion.

First of all, the quasi-symmetric algebra is functorial in the sense that 
any homomorphism of rings $R \to S$ compatible with epimorphisms to $A$ 
induces an epimorphism
$$\xymatrix{
\qSym_{R\to A}(\mathfrak a) \ar@{->>}[r] & \qSym_{S\to A}(\mathfrak a)
}\quad.$$

Indeed, the homomorphisms $R \to S \to A$ induce the 
middle row in the diagram
$$\xymatrix@=12pt{
K_R \ar[d] \ar[r] & K_S \ar[d] \\
\Sym_R(I) \ar[r] \ar[d] & \Sym_S(J) \ar[r] \ar[d] & \Sym_A(\mathfrak a) \\
\Rees_R(I) \ar[r] & \Rees_S(J)
}$$
where $J$ is the inverse image of $\mathfrak a$ in $S$, and $K_R$, $K_S$ 
are the kernels of the vertical maps to the Rees algebras. Since 
$K_R\Sym_A(\mathfrak a)\subset K_S\Sym_A(\mathfrak a)$, there is an induced 
epimorphism
$$\xymatrix{
\qSym_{R\to A}(\mathfrak a)=\frac{\Sym_A(\mathfrak a)}
{K_R\Sym_A(\mathfrak a)} \ar@{->>}[r] & 
\frac{\Sym_A(\mathfrak a)}{K_S\Sym_A(\mathfrak a)}
=\qSym_{S\to A}(\mathfrak a)
}\quad.$$

Pictorially, we have the commutative diagram:
$$\xymatrix{
\Sym_R(I) \ar[d] \ar[r] & \ar@{}[dr] |{\boxed{}}
\Sym_S(J) \ar[d] \ar[r] & \Sym_A(\mathfrak a) \ar[d] \\
\Rees_R(I) \ar[r] & \Rees_S(J) \ar[r] & \qSym_{S\to A}(\mathfrak a)
}$$
where the square on the right is cocartesian by definition. As
$\qSym_{R\to A}(\mathfrak a)$ satisfies a universal property (as a
tensor product) there is an induced canonical homomorphism
$\qSym_{R\to A}(\mathfrak a) \to \qSym_{S\to A}(\mathfrak a)$.

\subsection{}
The functoriality is the key to most of the following remarks, whose proof 
is left to the reader.

\begin{lemma}\label{standard}
Let $R\to A$, $\mathfrak a$, $I$ be as above.
\begin{enumerate}
\item The quasi-symmetric algebra bridges between the Rees algebra and the 
ordinary symmetric algebra of $\mathfrak a$ in $A$:
$$\xymatrix{
\Sym_A(\mathfrak a) \ar@{->>}[r] &
\qSym_{R\to A}(\mathfrak a) \ar@{->>}[r] & 
\Rees_A(\mathfrak a)
}\quad.$$
\item\label{ci} If $\Sym_R(I)=\Rees_R(I)$, then $\qSym_{R\to A}(\mathfrak a) 
=\Sym_A(\mathfrak a)$.
\item If $R \to A$ splits, then $\qSym_{R\to A}(\mathfrak a) 
=\Rees_A(\mathfrak a)$.
\item There is an epimorphism $\xymatrix{\Rees_R(I) \ar@{->>}[r] & 
\qSym_{R\to A}(\mathfrak a)}$.
\end{enumerate}
\end{lemma}

\begin{example} 
If $I$ is a complete intersection in $R$, then $\qSym_{R\to A}(\mathfrak a) 
=\Sym_A(\mathfrak a)$ by part~\ref{ci}.~in Lemma~\ref{standard} (since then 
the symmetric and Rees algebras of $I$ in $R$ coincide, \cite{MR31:1275}).

This shows that $\qSym_{R\to A}(\mathfrak a)$ may depend on $R$. For example, 
let $A=\Cbb[x,y]/(xy)$, $\mathfrak a=(x,y)$, $R=\Cbb[x,y]$; then
$$\qSym_{R\to A}(\mathfrak a)=\Sym_A(\mathfrak a)\ne\Rees_A(\mathfrak a)=
\qSym_{A\to A}(\mathfrak a)
\quad.$$
\end{example}

However, one of the main results of this section (Theorem~\ref{concom}) 
will show that $\qSym_{R\to A}(\mathfrak a)$ is in fact independent of $R$ 
provided that $R$ is constrained to be {\em regular\/.}

\subsection{}
There are two important cases in which the induced epimorphism is in fact 
an isomorphism.

\begin{lemma}\label{indepalg}
Let $R\to S$ be a ring homomorphism compatible with epimorphisms $R\to A$ 
and $S\to A$; let $\mathfrak a$ be an ideal of $A$, and let $I$, $J$ 
resp.~be the inverse images of $\mathfrak a$ in $R$, $S$. Then the induced 
epimorphism
$$\xymatrix{
\qSym_{R\to A}(\mathfrak a) \ar@{->>}[r] & \qSym_{S\to A}(\mathfrak a)
}$$
is an isomorphism if
\begin{enumerate}
\item the homomorphism $R\to S$ splits; or
\item $S$ is $R$-flat, and $J=IS$.
\end{enumerate}
\end{lemma}

\begin{proof}
In the first situation, if a composition $R \to S\to R$ is the 
identity we obtain a decomposition of the identity
$$\xymatrix{
\qSym_R(I) \ar@{->>}[r] & \qSym_S(J) \ar@{->>}[r] & \qSym_R(I)
}$$
implying that both maps are isomorphisms.

In the second situation, since $S$ is flat over $R$ we have $I^m S=S\otimes_R 
I^m$ for all $m$. Thus
$$\Rees_S(J)=\Rees_S(IS)=\oplus_m S\otimes_R I^m=S\otimes_R (\oplus_m I^m)
=S\otimes_R \Rees_R(I)\quad.$$
On the other hand, and again using flatness,
$$\Sym_S(J)=\Sym_S(IS)=\Sym_S(I\otimes_R S)=S\otimes_R \Sym_R(I)\quad,$$
by \cite{Bourbaki}, III \S6, Proposition~7. Thus
$$\aligned
\Sym_S(J)\otimes_{\Sym_R(I)} \Rees_R(I)
&= \big(S\otimes_R \Sym_R(I)\big) \otimes_{\Sym_R(I)} \Rees_R(I)\\
&= S\otimes_R \big(\Sym_R(I) \otimes_{\Sym_R(I)} \Rees_R(I)\big)\\
&= S\otimes_R \Rees_R(I)\\
&= \Rees_S(J)\quad.
\endaligned$$
This shows that the square on the left in the diagram at the end of 
\S\ref{definqsym} is cocartesian, implying the assertion.
\end{proof}

\subsection{}
We now move to the geometric setting. All our schemes are of finite type 
over a field $k$.

Let $Y\subset X\subset M$ be closed embeddings of schemes. We denote by
$\mathcal J_{Y,X}$, resp.~$\mathcal J_{Y,M}$ the ideals of $Y$ in $X$
and $M$, respectively.

\begin{defin}\label{qsymdef}
The {\em quasi-symmetric algebra\/} $\qSym_{X\subset M}(\mathcal
J_{Y,X})$ is the graded $\mathcal O_X$-algebra
$$\qSym_{X\subset M}(\mathcal J_{Y,X}):=\Sym_{\mathcal O_X}
(\mathcal J_{Y,X})\otimes_{\Sym_{\mathcal O_M}(\mathcal J_{Y,M})}
\Rees_{\mathcal O_M}(\mathcal J_{Y,M})\quad.$$
\end{defin}

In other words, $\qSym_{X\subset M}(\mathcal J_{Y,X})$ sheafifies the 
local construction given by Definition~\ref{qsym}. Every commutative 
diagram
$$\xymatrix{
& & N \ar@{->>}[d]^\pi \\
Y \ar@{^{(}->}[r] & X \ar@{^{(}->}[ru]^j \ar@{^{(}->}[r]^i & M
}$$
determines an epimorphism
$$\xymatrix{
\qSym_{X\subset M}(\mathcal J_{Y,X}) \ar@{->>}[r] &
\qSym_{X\subset N}(\mathcal J_{Y,X})
}$$
and we are interested in conditions guaranteeing that this map is an 
isomorphism.

\begin{lemma}\label{indepgeom}
The induced epimorphism is an isomorphism if
\begin{enumerate}
\item $N=M\times \Abb^n$; or
\item $\pi$ is flat, and $j(X)$ is a connected component of $\pi^{-1}(i(X))$.
\end{enumerate}
\end{lemma}

\begin{proof}
These follow from Lemma~\ref{indepalg}. As the matter can be checked 
locally, we may assume $M=\Spec R$, $N=\Spec S$, $X=\Spec A$, $Y$ is given 
by an ideal $\mathfrak a$ in $A$, and we have a commutative diagram
$$\xymatrix{
& & S \ar@{->>}[ld] \\
A/\mathfrak a & A \ar@{->>}[l] & R \ar@{->>}[l] \ar@{^{(}->}[u] 
}$$
Denote by $K$, $L$ resp.~the kernels of $R\to A$, $S\to A$ resp.; and by 
$I$, $J$ resp.~the inverse images of $\mathfrak a$ in $R$, $S$ resp.

In the first situation $S=R[u_1,\dots,u_s]$ is a polynomial ring, and the 
splitting needed in order to apply Lemma~\ref{indepalg} holds because if 
$K$ is an ideal of $R$ then any left-inverse of the inclusion 
$R/K\hookrightarrow R/K[u_1,\dots,u_s]$ lifts to a left-inverse of 
$R\hookrightarrow R[u_1,\dots,u_s]$.

In the second situation, by hypothesis $S$ is flat over $R$, and there exists 
an $f\in S$ such that the epimorphism $S \to A$ lifts to an epimorphism $S_f 
\to A$ from the localization of $S$ at $f$, with kernel $KS_f=LS_f$. A 
fortiori $IS_f=JS_f$ is the inverse image of $\mathfrak a$ in $S_f$. As $S_f$ 
is flat over both $S$ and $R$, two applications of part~2.~from 
Lemma~\ref{indepalg} give the assertion.
\end{proof}

\begin{theorem}\label{smooth}
If $\pi:N\to M$ is a smooth map compatible with closed embeddings $X\subset 
M$, $X\subset N$, then for all closed subschemes $Y\subset X$ the induced 
epimorphism
$$\xymatrix{
\qSym_{X\subset M}(\mathcal J_{Y,X}) \ar@{->>}[r] &
\qSym_{X\subset N}(\mathcal J_{Y,X})
}$$
is an isomorphism.
\end{theorem}

\begin{proof}
Again the matter can be checked locally, so as $\pi$ is smooth we may assume 
that it can be written as a composition
$$\xymatrix{
N \ar[r]^-{\text{\'etale}} & M\times \Abb^s \ar[r] & M
}\quad;$$
by Lemma~\ref{indepgeom}, part~1., we may assume that $\pi$ itself is \'etale. 
In this case $\pi^{-1}(X) \to X$ is an \'etale map with a section; hence 
the image of $X$ in $N$ must be a connected component of $\pi^{-1}(X)$. 
As \'etale maps are flat, part~2.~in Lemma~\ref{indepgeom} concludes the 
proof.
\end{proof}

Theorem~\ref{smooth} shows that the quasi-symmetric algebras of $X$
collect into classes detecting specific `qualities' of the embeddings
$X\subset M$.  For example, if $X\subset M$ is a section of a smooth
projection $M\to X$ then $\qSym_{X\subset M}(\mathcal
J_{Y,X})=\Rees_{\mathcal O_X}(\mathcal J_{Y,X})$ for all closed
subschemes $Y\subset X$. In fact, only the features of the embedding
$X\subset M$ {\em near~$Y$\/} affect the corresponding quasi-symmetric
algebra.

\subsection{}
It is time to remove the dependence on the choice of an embedding 
$X\subset M$. For given $Y\subset X$, the epimorphisms on quasi-symmetric 
algebras induced by concatenation of embeddings $X\subset M\subset N$ 
make $\{\qSym_{X\subset M}(\mathcal J_{Y,X})\}_M$ into an inverse system.

\begin{defin}\label{absqsym}
Let $Y\subset X$ be a closed embedding of schemes. The {\em quasi-symmetric 
algebra\/} of $\mathcal J_{Y,X}$ is defined as the inverse limit
$$\qSym_{\mathcal O_X}(\mathcal J_{Y,X}):=\varprojlim_{M\supset X} 
\qSym_{X\subset M}(\mathcal J_{Y,X}) \quad.$$

The {\em quasi-symmetric blow-up\/} of $X$ along $Y$ is defined as the 
Proj of the quasi-symmetric algebra:
$$\qBl_YX:=\Proj(\qSym_{\mathcal O_X}(\mathcal J_{Y,X}))\quad.$$
\end{defin}

The quasi-symmetric blow-up carries a tautological line bundle $\mathcal 
O(-1)$, as do the more conventional $sB\ell_YX=\Proj(\Sym_{\mathcal O_X} 
(\mathcal J_{Y,X}))$ and $B\ell_YX=\Proj(\Rees_{\mathcal O_X}(\mathcal 
J_{Y,X}))$. Also, note that by Lemma~\ref{standard}, part~1., there are 
closed embeddings
$$B\ell_YX \subset \qBl_YX \subset sB\ell_YX\quad.$$

Theorem~\ref{smooth} is the key to the following concrete computation of the 
`absolute' quasi-symmetric algebra and blow-up.

\begin{theorem}\label{concom}
Let $Y\subset X\subset M$ be closed embeddings of schemes, with $M$ 
nonsingular. Then the canonical epimorphism
$$\xymatrix{
\qSym_{\mathcal O_X}(\mathcal J_{Y,X}) \ar@{->>}[r] &
\qSym_{X\subset M}(\mathcal J_{Y,X})
}$$
is an isomorphism.
\end{theorem}

\begin{proof}
The matter is local. Since locally every scheme is embedded in a 
nonsingular variety, it suffices to show that if $X\subset M\subset N$ are 
closed embeddings, with $M$ and $N$ nonsingular varieties, then 
$\qSym_{X\subset N}(\mathcal J_{Y,X})\to \qSym_{X\subset M} (\mathcal 
J_{Y,X})$ is an isomorphism. Factoring the embedding $M\subset N$ through 
the product, we have the diagram
$$\xymatrix@C=3pt{
& & & & & M\times N \ar@{->>}[dr] \\
X \ar@{^{(}->}[urrrrr] \ar@{^{(}->}[rrrr] &
& & & M \ar@{^{(}->}[ur] \ar@{^{(}->}[rr] & & N
}$$
which induces the commutative diagram of $\qSym$ algebras
$$\xymatrix@C=0pt{
& \qSym_{X\subset M\times N}(\mathcal J_{Y,X}) \ar@{->>}[ld]_\sim \\
\qSym_{X\subset M}(\mathcal J_{Y,X}) & &
\qSym_{X\subset N}(\mathcal J_{Y,X}) \ar@{->>}[ll] \ar@{->>}[ul]_\sim 
}$$ 
The diagonal arrow on the left is an isomophism because the diagonal 
embedding splits; the diagonal arrow on the right is an isomorphism by 
Theorem~\ref{smooth}. Thus the horizontal arrow is an isomorphism, as 
needed.
\end{proof}

\subsection{}
By Theorem~\ref{concom}, the inverse system of algebras $\qSym_{X\subset M}
(\mathcal J_{Y,X})$ stabilizes at nonsingular ambient varieties $M$. In 
fact, by part~4.~in Lemma~\ref{standard} there is a canonical embedding
$$\qBl_Y X\subset B\ell_Y M$$
induced by the surjection $\Rees_{\mathcal O_M}(\mathcal J_{Y,M}) \to 
\qSym_{\mathcal O_X}(\mathcal J_{Y,X})$; the line bundle $\mathcal O(-1)$ 
is the restriction of the line bundle of the exceptional divisor. 

Theorem~\ref{concom} implies immediately that $\qBl_YX$ fulfills the 
promise made in the introduction.

\begin{corol}\label{prom}
The quasi-symmetric blow-up $\qBl_YX$ is the largest subscheme of $sB\ell_YX$ 
which admits an embedding in $B\ell_YM$ (compatibly with the projection 
to $X$) for some scheme $M$ containing $X$.
\end{corol}

\begin{proof}
By Theorem~\ref{concom}, the quasi-symmetric blow-up $\qBl_YX$ can be embedded 
in $B\ell_YM$ for any {\em nonsingular\/} variety $M$ containing $X$.

On the other hand, if $\mathcal S$ is any quotient algebra of 
$\Sym_{\mathcal O_X}(\mathcal J_{Y,X})$ which is also a quotient of 
$\Rees_{\mathcal O_M}(\mathcal J_{Y,M})$ (for some $M$) then there is an 
induced epimorphism $\qSym_{X\subset M}(\mathcal J_{Y,X})\to \mathcal S$ 
by definition of $\qSym$ (as a tensor product). Hence we obtain a surjection 
$\qSym_{\mathcal O_X}(\mathcal J_{Y,X})\to \mathcal S$, showing that 
$\Proj(\mathcal S)\subset \qBl_YX$, as needed.
\end{proof}

It is natural to ask whether the embedding $\qBl_YX\subset B\ell_YM$ (for 
$M$ a nonsingular variety containing $X$) can be realized concretely, just as 
the embedding $B\ell_Y X\subset B\ell_Y M$ can be realized as a `proper 
transform'.

\begin{defin}
Let $Y\subset X\subset M$ be closed embeddings of schemes. The {\em 
principal transform\/} of $X$ in the blow-up $B\ell_Y M \overset\rho\to M$ of 
$M$ along $Y$ is the residual to the exceptional divisor in $\rho^{-1}(X)$ 
(in the sense of \cite{MR85k:14004}, Definition~9.2.1).
\end{defin}

This definition would appear to depend on $M$; at any rate, $\rho^{-1}(X)$ 
certainly depends on $M$ as it contains the exceptional divisor of 
$B\ell_YM$. However, the next result claims that the {\em principal\/} 
transform is almost as intrinsic to $X$, $Y$ as is the {\em proper\/} 
transform.

\begin{theorem}\label{printran}
Let $Y\subset X\subset M$ be closed embeddings of schemes, with $M$ a 
nonsingular variety. Then the quasi-symmetric blow-up of $X$ along $Y$ 
equals the principal transform of $X$ in $B\ell_YM$.
\end{theorem}

This is an easy consequence of the following computation of $\qSym_{\mathcal 
O_X}(\mathcal J_{Y,X})$.

\begin{lemma}
With notations as in the statement of the theorem, 
$$\qSym_{\mathcal O_X}(\mathcal J_{Y,X})=\oplus_{d\ge 0} 
\frac{\mathcal J_{Y,M}^d}{\mathcal J_{X,M}\mathcal J_{Y,M}^{d-1}}$$
(where we set $\mathcal J_{Y,M}^{-1}=\mathcal O_M$).
\end{lemma}

\begin{proof}
By Theorem~\ref{concom} we have $\qSym_{\mathcal O_X}(\mathcal J_{Y,X}) 
\cong\qSym_{X\subset M}(\mathcal J_{Y,X})$; we compute the latter.
For $d\ge 1$ set up the commutative diagram with exact rows
$$\xymatrix{
&  &  & 0 \ar[d]\\
&  & {\Tors_d} \ar[r] \ar[d] & {\Disc_d} \ar[r] \ar[d] & 0\\
& {J_{X,M}\cdot \Sym^{d-1}\mathcal J_{Y,M}} \ar[r] \ar[d] & 
{\Sym^d\mathcal J_{Y,M}} \ar[r]
\ar[d] & {\Sym^d({\mathcal J_{Y,M}}/\mathcal J_{X,M})} \ar[r] \ar[d] & 0\\
0 \ar[r] & {\mathcal J_{X,M}\cdot \mathcal J_{Y,M}^{d-1}} \ar[r]
\ar[d] & {\mathcal J_{Y,M}^d} \ar[r]
\ar[d] & {\mathcal J_{Y,M}^d/(\mathcal J_{X,M}\cdot \mathcal J_{Y,M}^{d-1})}
\ar[r] \ar[d] & 0\\
& 0 & 0 & 0
}$$
where all $\Sym$ are over $\mathcal O_M$, $J_{X,M}$ denotes the
image of $\mathcal J_{X,M}$ in $\Sym^1\mathcal J_{Y,M}=\mathcal
J_{Y,M}$, $\Tors_d$ is defined to make the central column exact, $\Disc_d$ 
is its image in $\Sym^d({\mathcal J_{Y,M}}/\mathcal J_{X,M})$. A diagram 
chase shows that the column on the right is exact. This gives
$$\mathcal O_M\oplus\oplus_{d\ge 1} \frac{\mathcal J_{Y,M}^d}{\mathcal J_{X,M}
\mathcal J_{Y,M}^{d-1}}=\Sym_{\mathcal O_M}({\mathcal J_{Y,M}}/\mathcal J_{X,M}) 
\otimes_{\Sym_{\mathcal O_M} \mathcal J_{Y,M}} \Rees_{\mathcal O_M}(\mathcal 
J_{Y,M})\quad.$$
Tensoring by $\mathcal O_X$ only affects the term of degree~0 on the left. 
On the other hand,
$$\mathcal O_X\otimes_{\mathcal O_M}\Sym_{\mathcal O_M}(\mathcal J_{Y,M}/
\mathcal J_{X,M})=\Sym_{\mathcal O_X}(\mathcal J_{Y,X})\quad;$$
thus
$$\oplus_d \frac{\mathcal J_{Y,M}^d}{\mathcal J_{X,M}\mathcal 
J_{Y,M}^{d-1}}=\qSym_{\mathcal O_X}(\mathcal J_{Y,X})$$
by the associativity of tensor products.
\end{proof}

\subsection{}
For $Y\subset X\subset M$, and $M$ not necessarily nonsingular, we can of 
course consider a quasi-symmetric blow-up $\Proj(\qSym_{X\subset 
M}(\mathcal J_{Y,X}))$. The analog of Theorem~\ref{printran} realizes this 
blow-up as the principal transform of $X$ in the (ordinary) blow-up of $M$ 
along $Y$. 

This more general notion will not be used in the rest of this paper.
We observe that every such blow-up is contained in the quasi-symmetric
blow-up of Definition~\ref{absqsym}, since any $M$ is contained
locally in a nonsingular variety. The ordinary blow-up $B\ell_YX$ is
recovered for $M=$ (for example) $X\times \Abb^1$.

The example of $X=$ three noncoplanar lines through a point $p=Y$ in 
projective space shows that $qB\ell_YX$ may have components of higher 
dimension than $X$. However, it is easily checked that if the lines are 
coplanar then the quasi-symmetric blow-up has dimension~1 (it consists of 
three disjoint lines union a double line connecting them). More generally:

\begin{corol}
If $X$ can locally be embedded as a hypersurface in a nonsingular
irreducible variety, then for every $Y\subset X$ the quasi-symmetric
blow-up $\qBl_YX$ is equidimensional.
\end{corol}

\begin{proof}
Immediate consequence of Theorem~\ref{printran}.
\end{proof}

Hypersurfaces of nonsingular varieties will be our main concern in the 
rest of the paper.


\section{The conormal and characteristic cycles of a hypersurface}
\label{cycles}

\subsection{} 
We now move from the generalities in \S\ref{qsymblowups} to our 
application to the theory of Chern classes of singular varieties. In this 
section we will deal with the theory at the level of Lagrangian cycles in 
the cotangent bundle of an ambient nonsingular variety; in the next 
section we will extract the information more closely pertaining to 
characteristic classes.

Our main objective in this section is to show that the notion
introduced in \S\ref{qsymblowups} gives a concrete realization of the
characteristic cycle of a hypersurface $X$ in a nonsingular ambient
variety $M$. In a nutshell, the {\em characteristic\/} cycle of $X$ is
the cycle of the {\em quasi-symmetric blow-up\/} of $X$ along its
singularity subscheme. This fact should be appreciated in conjunction
with the (straightforward) observation that the {\em conormal\/} cycle
of $X$ is the cycle of its ordinary blow-up along the same subscheme.

Realizing the characteristic cycle allows us to give a direct computation 
of the Chern-Schwartz-MacPherson classes of a hypersurface, following the 
same philosophy behind other characteristic classes (specifically the 
classes introduced in \cite{MR82c:14017} and \cite{MR85k:14004}, 
Example~4.2.6). This requires a certain care in handling the appropriate 
tautological line bundles; we work this out in \S\ref{charcla}.

After the preliminary work done in \S\ref{qsymblowups}, the main result in 
this section follows easily from the existing literature on characteristic 
classes for singular hypersurfaces.

In this section we also identify a condition under which the 
quasi-symmetric blow-up needed here equals the symmetric blow-up. In this 
situation, the Chern-Schwartz-MacPherson class of the hypersurface can be 
efficiently expressed in terms of the Chern class of a certain coherent 
sheaf defined on it.

\subsection{} We work over an algebraically closed field of
characteristic~0. Throughout the rest of the paper $M$ will denote a
nonsingular irreducible algebraic variety, and $X$ will be the
zero-scheme of a nonzero section $F$ of a line bundle $\mathcal L$ on
$M$; we will say that $X$ is a {\em hypersurface\/} for short. For
convenience we will implicitly assume that $X$ is {\em reduced,\/}
although this is not an essential requirement (cf.~\S\ref{nonred}).

The singularity locus of $X$ has an interesting, possibly nonreduced
scheme structure. We will denote by $Y$ this {\em singularity
subscheme\/} of $X$ (see \S\ref{conorm} for the precise definition).

We begin by recalling several well-established notions, for the benefit of 
the non-expert and in order to establish notations. The informed and 
impatient reader can safely skip to \S\ref{skiphere}.

\subsection{}\label{gencon}
A {\em constructible function\/} on a variety $V$ is a finite linear
combination
$$\sum n_W\,\one_W$$
where the summation ranges over (closed, irreducible) subvarieties
$W\subset V$, $n_W\in\Zbb$, and $\one_W$ denotes the function that is
the constant $1$ on $W$, and $0$ outside of $W$. We denote by $C(V)$
the group of constructible functions on $V$. If $f:V_1\to V_2$ is a
proper map, a push-forward $C(f):C(V_1)\to C(V_2)$ is defined by
setting, for $W$ a subvariety of $V_1$ and $p\in V_2$,
$$C(f)(\one_W)(p)=\chi(f^{-1}(p)\cap W)\quad,$$
and extending by linearity. Here $\chi$ denotes the topological Euler
characteristic when working over $\Cbb$; see \cite{MR91h:14010}, \S3,
for the extension of the theory to arbitrary algebraically closed
field of characteristic~0.

With this push-forward, the assignment
$$\mathcal C\quad:\quad V \mapsto C(V)$$
yields a covariant functor from algebraic varieties to abelian groups. 

\subsection{} A fundamental result of Robert MacPherson
(\cite{MR50:13587} and \cite{MR91h:14010}) compares this functor to
the functor 
$$\mathcal A\quad:\quad V \mapsto A(V)$$
assigning to a variety its Chow group: there exists a natural
transformation
$$c_*\quad:\quad \mathcal C \leadsto \mathcal A$$
such that, for $V$ a {\em nonsingular\/} variety, the induced group
homomorphism
$$C(V) \rightarrow A(V)$$
maps $\one_V$ to the total Chern class of the tangent bundle of $V$:
$$\one_V \mapsto c(TV)\cap [V]\quad.$$
For arbitrarily singular $V$, one may then define a (total) `Chern class'
in the Chow group of $V$, by setting
$$\csm(V):= c_*(\one_V)\quad;$$
thus $\csm(V)=c(TV)\cap[V]$ if $V$ is nonsingular.
Jean-Paul Brasselet and Marie-H\'el\`ene Schwartz later discovered
that this class defined by MacPherson is in fact Alexander dual to a
class previously defined by Schwartz (\cite{MR32:1727},
\cite{MR35:3707}; and \cite{MR83h:32011}); nowadays, $\csm(V)$ is
commonly named the {\em Chern-Schwartz-MacPherson class\/} of $V$.

\subsection{}\label{sabbah} 
In MacPherson's approach, the natural transformation $c_*$ is defined 
directly by requiring that 
$$c_*(\Eu_V)=\cma(V)$$
for all varieties $V$. Here $\cma(V)$ stands for the {\em
Chern-Mather\/} class of $V$, and $\Eu_V$ is the {\em local Euler
obstruction,\/} a measure of the singularities of $V$. Both these
notions were defined in \cite{MR50:13587}, and have been the
subject of intense study since; again we refer the reader to
\cite{MR91h:14010} for a very readable treatment and for the extension
of the theory to arbitrary algebraically closed fields of
characteristic~0.

A different approach to the definition of $c_*$ emerged in the work of
Claude Sabbah (\cite{MR87f:32031}, \cite{MR91h:14010}, and
\cite{MR1795550}, \S1). The natural transformation $c_*$ can be
obtained (up to taking a harmless dual) as the composition of two 
transformations
$$\mathcal C \leadsto \mathcal L\leadsto \mathcal A\quad,$$
where $\mathcal L$ denotes the functor assigning to a variety $V$ the
group $L(V)$ of {\em Lagrangian cycles\/} over $V$, with a suitably
defined push-forward. If $V\subset M$ is an embedding of $V$ into a
nonsingular variety $M$, the Lagrangian cycles over $V$ are the
Lagrangian cycles in the restriction $\Pbb(T^*M)|_V$ of the
projectivized cotangent bundle of $M$. As is well known
(\cite{MR91h:14010}, Lemma~3), every Lagrangian subvariety
over $V$ is in fact the projective conormal space $\Pbb(T^*_WM)$ of a
closed subvariety $W\subset V$. Hence $L(V)$ is the free abelian group
on the set of subvarieties of $V$; the realization of $L(V)$ as a
group of cycles in $\Pbb(T^*M)|_V$, for some nonsingular $M$
containing $V$, yields a good notion of push-forward of elements of
$L(V)$ (see p.~2829-31 in \cite{MR91h:14010} for details). 

The second step $\mathcal L\leadsto \mathcal A$ in the above
decomposition can be expressed in terms of standard intersection
theory, and will be recalled in \S\ref{sabbah2}. The first step,
$\mathcal C\leadsto \mathcal L$, is considerably subtler. It is
determined by the requirement that, for all (closed, irreducible)
subvarieties $W\subset V$, the local Euler obstruction of $W$
correspond (up to a sign) to the {\em conormal cycle\/} of $W$ in $M$:
$$(-1)^{\dim W}\Eu_W \mapsto [\Pbb(T^*_WM)]$$
For every constructible function $\varphi\in C(V)$ we obtain then a
{\em characteristic cycle\/}
$$\Ch(\varphi)\in L(V)\quad.$$

The cycle $\Ch(\one_V)$ (realized as above, that is, in terms of an
embedding $V\subset M$) is called the {\em characteristic cycle of
$V$ (in $M$).\/}

\subsection{}\label{skiphere}
Summarizing, there are two important cycles associated to a variety $V$ in 
the projectivized cotangent bundle $\Pbb(T^*M)=\Proj(\Sym_M((\Omega^1_M)^\vee))$ 
of any nonsingular variety $M$ in which $V$ is embedded:\begin{itemize}

\item the conormal cycle $[\Pbb(T_V^*M)]$, corresponding (up to sign)
to the local Euler obstruction of $V$, and to the Chern-Mather class
of $V$; and

\item the characteristic cycle $\Ch(V)$ of $V$, likewise corresponding
to the constant function $\one_V$ and to the Chern-Schwartz-MacPherson
class of $V$.\end{itemize}

Explicitly realizing $\Ch(V)$ `from the definition' requires finding 
subvarieties $W$ of $V$ and integers $e_W$ such that $\one_V=\sum_W e_W 
\Eu_W$. This information is extremely subtle. `Index formulas' 
(cf.~\cite{MR83a:32010}) provide an approach to extracting this
information, but we do not know of any computationally effective
method to implement such formulas.

Our goal here is the construction of a scheme whose cycle is the 
characteristic cycle of a hypersurface $X$ of a nonsingular variety $M$. In 
principle this construction can be performed by symbolic computation programs 
such as {\tt Macaulay2.} An entirely analogous realization of the conormal 
cycle is more readily available, and will be recalled in a moment.

The theory recalled above applies to varieties, and in particular 
requires $V$ to be {\em reduced.\/} Because of this, we will assume that 
our hypersurfaces are reduced in what follows (but see \ref{nonred} below).

\subsection{}\label{properplace} 
According to the framework recalled above, the conormal and characteristic 
cycles arise as cycles in the (projectivized) cotangent bundle of an 
ambient nonsingular variety. It is our opinion that these cycles have a 
right to exist freely, independent of an ambient variety; but we will wait 
until \S\ref{charcla} to fully make this point. For the time being we will 
house the cycles in the usual place, which amounts to finding an 
appropriate ambient for the blow-ups considered in \S\ref{qsymblowups}.

The section $F$ of $\mathcal L$ defining $X$ determines a section $s$ 
of the bundle $\mathcal P^1_M\mathcal L$ of principal parts of 
$\mathcal L$:
$$s:\mathcal O_M \rightarrow \mathcal P^1_M\mathcal L\quad;$$
we let $Y$ denote the zero-scheme of $s$ in $M$, and we call $Y$ the
`singularity subscheme' of $X$. Composing $s$ with the projection to
$\mathcal L$ recovers $F$:
$$\xymatrix{
{\mathcal O_M} \ar[r]_s \ar@/^1pc/[rr]^F & {\mathcal P^1_M\mathcal L}
\ar[r] & {\mathcal L}}
\quad;$$
hence $s$ induces a section of $\Omega^1_M\otimes\mathcal L$ on $X$,
which is natural to name $dF$:
$$dF:\mathcal O_X\longrightarrow (\Omega^1_M\otimes\mathcal L)|_X
\quad;$$
the subscheme $Y$ is the zero-scheme of $dF$ on $X$. It is easily
checked that, locally, $dF$ is given by the partial derivatives of $F$
with respect to a set of local parameters for $M$; hence $Y$ is
supported on the singular locus of $X$, justifying its name. Locally,
we can write (abusing notations):
$$\mathcal J_{Y,M}=\left(F,\frac{\partial F}{\partial x_1}, \dots,
\frac{\partial F}{\partial x_n}\right)$$
for the ideal or $Y$ in $M$. We will write $(F)$ for the ideal of $X$ in 
$M$, as this is given by the vanishing of the section $F$ of $\mathcal L$.

Dualizing $s: \mathcal O_M \to \mathcal P^1_M\mathcal L$ we get an 
epimorphism
$$\xymatrix{
(\mathcal P^1_M\mathcal L)^\vee \ar@{->>}[r] & \mathcal J_{Y,M}
}$$
and from this, Lemma~\ref{standard}, and Theorem~\ref{concom} 
the epimorphisms
$$\xymatrix{
\Sym_{\mathcal O_M}((\mathcal P^1_M\mathcal L)^\vee) \ar@{->>}[r] &
\Rees_{\mathcal O_M}(\mathcal J_{Y,M}) \ar@{->>}[r] &
\qSym_{\mathcal O_X}(\mathcal J_{Y,X})
}\quad.$$
Since $\qSym_{\mathcal O_X}(\mathcal J_{Y,X})$ is an $\mathcal 
O_X$-module, tensoring by $\mathcal O_X$ gives an epimorphism
$$\xymatrix{
\Sym_{\mathcal O_M}((\mathcal P^1_M\mathcal L)^\vee|_X) \ar@{->>}[r] &
\qSym_{\mathcal O_X}(\mathcal J_{Y,X})
}\quad.$$
Finally, composing with $\mathcal L^\vee \hookrightarrow (\mathcal
P^1_M\mathcal L)^\vee$ gives the zero-map over $X$, showing that there
is a surjection $$\xymatrix{
\Sym((\Omega^1_M\otimes \mathcal L)^\vee|_X) \ar@{->>}[r] &
\qSym_{\mathcal O_X}(\mathcal J_{Y,X})
}\quad.$$

Since $\qSym_{\mathcal O_X}(\mathcal J_{Y,X})$ dominates all 
quasi-symmetric algebras of $\mathcal J_{Y,X}$, and in particular the Rees 
algebra, this shows (taking Proj) that there are closed embeddings
$$B\ell_YX \subset \qBl_YX \subset \Proj(\Sym((\Omega^1_M\otimes \mathcal 
L)^\vee|_X))=\Pbb(T^*M\otimes\mathcal L|_X)\cong \Pbb(T^*M|_X)\quad.$$

\subsection{}\label{conorm} 
The following statement is only one step away from the definitions, but it is 
excellent preparation for the main result of the section, Theorem~\ref{qSym} 
below.

\begin{theorem}\label{rees}
The conormal cycle $[\Pbb(T_X^*M)]$ of $X$ in $M$ equals
$$[B\ell_YX]=[\Proj(\Sym_{X\subset X}(\mathcal J_{Y,X}))]\quad.$$
\end{theorem}

\begin{proof} Recall that we are assuming that $X$ is reduced. The
conormal space $T^*_XM$ of $X$ in $M$ is the closure in $T^*M$ of the
kernels of the projection
$$(T^*M)_x \longrightarrow (T^*X)_x \longrightarrow 0$$
over nonsingular points $x$ of $X$. In other words, the projectivized
conormal space of $X$ is the closure of the image of the section
$$X^{\text{reg}} \longrightarrow
\Pbb(T^*M)|_X=\Pbb(T^*M\otimes\mathcal L)|_X$$
induced on the set $X^{\text{reg}}$ of regular points of $X$ by the
section $dF$ determined above. Chasing the morphisms collected above shows 
that this is precisely how $B\ell_YX$ is embedded in $\Pbb(T^*M|_X)$ over 
regular points of $X$. Hence $B\ell_YX$ and the projectivized conormal 
space agree over regular points of $X$, and it follows that they agree 
everywhere, as needed.
\end{proof}

\subsection{}\label{ss:qSym}
The next result is our main application of the construction developed in 
\S\ref{qsymblowups}; it does for the characteristic cycle precisely what 
Theorem~\ref{rees} does for the conormal cycle.

\begin{theorem}\label{qSym}
The characteristic cycle $[\Ch(X)]$ of $X$ in $M$ equals
$$(-1)^{\dim X}[\qBl_YX]=(-1)^{\dim X}[\Proj(\qSym_{X\subset M}(\mathcal 
J_{Y,X}))]\quad.$$
\end{theorem}

The annoying sign is due to established (thus unavoidable) conventions, and 
reflects the fact that the lagrangian point of view is best suited to build 
a {\em co\/}tangent theory of characteristic classes.

Modulo the work done in \S\ref{qsymblowups}, the statement is an easy 
consequence of results in the literature on characteristic classes for 
singular varieties.

\begin{proof}
By Theorem~\ref{printran}, $[\qBl_YM]$ equals the principal transform of 
$X$ in $B\ell_YM$, so the claim is that the latter computes $\Ch(X)$, with 
due attention to the sign. Over $\Cbb$, this statement is Corollary~2.4 in 
\cite{MR1795550}; for arbitrary algebraically closed fields of 
characteristic~0, it can be obtained from Claim~2.1 in \cite{MR1819626}.
\end{proof}

\subsection{}\label{xcond}
We will now identify a technical condition under which the algebra 
$\qSym_X(Y)$ is nothing but the symmetric algebra of $\mathcal J_{Y,X}$. As 
a consequence of Theorem~\ref{qSym}, the characteristic cycle of hypersurfaces 
satisfying this condition is (up to sign) the cycle of the {\em symmetric\/} 
blow-up of their singularity subschemes. This both simplifies matters 
computationally (since packages such as {\tt Macaulay2\/} have built-in 
functions computing symmetric algebras) and is philosphically intriguing: 
in this case, the characteristic cycle is realized as the `linear fiber 
space' ({\em Linear Faserraum\/}, cf.~\cite{MR36:4024}) corresponding to the 
ideal sheaf $\mathcal J_{Y,X}$. While the fibers of the characteristic 
cycle are always linear, we do not know if every characteristic cycle can 
be realized as a linear fiber space.

As above, $F$ denotes the section of the line bundle $\mathcal L$ 
on $M$ whose zero-scheme is the hypersurface $X$.
For the purpose of this discussion, a {\em homogeneous, degree $d$
differential operator satisfied by $F$\/} is a local section of
$\Sym^{d}(\mathcal P^1_M\mathcal L)^\vee$ mapping to~$0$ in $\mathcal
J^d_{Y,M}$ via the map induced by the surjection $(\mathcal
P^1_M\mathcal L)^\vee \to \mathcal J_{Y,M}$ whose existence we pointed
out in \S\ref{properplace}. In terms of local parameters $x_1,\dots,x_n$
on $M$, this object is nothing but a homogeneous polynomial
$$P(T_0,\dots,T_n)$$
with coefficients (local) functions on $M$, such that
$$P\left(F,\frac{\partial F}{\partial x_1},\dots,\frac{\partial F}{\partial 
x_n}\right)\equiv 0\quad;$$
we will express the condition in this slightly imprecise but more vivid 
language, leaving to the reader the task of translating it into a global, 
coordinate-free formulation.

The simplest way to manufacture homogeneous differential operators of 
degree $d$ satisfied by $F$ is as a product
$$P=P_0\cdot T_0+P_1\cdot T_1+\dots+P_n\cdot T_n\quad,$$
where the $P_i$ are homogeneous polynomials of degree $d-1$ in
$T_0,\dots,T_n$, and 
$$P_0\cdot F+P_1\cdot \frac{\partial F}{\partial x_1}+\dots+P_1\cdot 
\frac{\partial F}{\partial x_n}=0\quad.$$
We say that such operators are {\em trivially\/} satisfied by $F$.
The {\em \x-condition\/} on $X$ is a softening of this requirement, on 
operators of sufficiently high degree satisfied by~$F$.

\begin{defin}
A hypersurface $X$ in $M$ {\em satisfies the \x-condition\/} if there 
exists a $d_0$ such that every homogeneous differential operator $P$ of 
degree $d\ge d_0$ and satisfied by $F$ can be written as
$$P=P_0\cdot T_0+P_1\cdot T_1+\dots+P_n\cdot T_n\quad,$$
where the $P_i$ are homogeneous polynomials of degree $d-1$ in
$T_0,\dots,T_n$, and 
$$P_1\cdot \frac{\partial F}{\partial x_1}+\dots+P_n\cdot 
\frac{\partial F}{\partial x_n}\in (F)\quad.$$
\end{defin}

\begin{prop}\label{symqsym}
A hypersurface $X$ satisfies the \x-condition if and only if
$$\Sym^d_X(Y)\cong \qSym^d_X(Y)$$
for $d\gg 0$.
\end{prop}

\begin{proof}
In the hypersurface case, we can complete the diagram in the proof of 
Theorem~\ref{printran} so that all rows and columns are exact:
$$\xymatrix@=14pt{
& 0 \ar[d] & 0 \ar[d] & 0 \ar[d]\\
0 \ar[r] & {T_0\cdot \Tors_{d-1}} \ar[r] \ar[d] & {\Tors_d} \ar[r]
\ar[d] & {\Disc_d} \ar[r] \ar[d] & 0\\
0 \ar[r] & {T_0\cdot \Sym^{d-1}\mathcal J_{Y,M}} \ar[r] \ar[d] & 
{\Sym^d\mathcal J_{Y,M}} \ar[r]
\ar[d] & {\Sym^d({\mathcal J_{Y,M}}/(F))} \ar[r] \ar[d] & 0\\
0 \ar[r] & {F\cdot \mathcal J_{Y,M}^{d-1}} \ar[r] \ar[d] & {\mathcal 
J_{Y,M}^d} \ar[r]
\ar[d] & {\mathcal J_{Y,M}^d/(F) \mathcal J_{Y,M}^{d-1}} \ar[r] \ar[d] 
& 0\\
& 0 & 0 & 0
}$$
(the leftmost column is exact as $F$ is a non-zero-divisor, and it follows 
that the top row is exact). We have to verify that $\Disc_d=0$ for $d\gg 0$ 
if and only if $X$ satisfies the \x-condition.

Now (cf.~for example \cite{MR95g:13005}, Chapter~2) $\Tors_d$ can be described 
as the space of degree-$d$ homogeneous operators satisfied by $F$, modulo 
those trivially satisfied by $F$. Hence 
$$\Disc_d=\frac{\Tors_d}{T_0\cdot \Tors_{d-1}}$$ 
is $0$ if and only if every degree-$d$ homogeneous operator satisfied by 
$F$ is equivalent to a multiple of $T_0$ modulo trivial ones.
That is, if and only if for every homogeneous $P$ of degree $d$ such that
$$P\left(F,\frac{\partial F}{\partial x_1},\dots,\frac{\partial F}{\partial 
x_n}\right)\equiv 0$$
there exists a $Q$, homogeneous of degree $d-1$ and such that
$$P-T_0\cdot Q=P_0\cdot T_0+P_1\cdot T_1+\dots+P_n\cdot T_n$$
with
$$P_0\cdot F+P_1\cdot \frac{\partial F}{\partial x_1}+\dots+
P_n\cdot \frac{\partial F}{\partial x_n}=0\quad.$$
It is straightforward to verify that this latter condition is satisfied 
for $d\gg 0$ if and only if $X$ satisfies the \x-condition.
\end{proof}

\begin{corol}\label{sym}
If $X$ satisfies the \x-condition, then the characteristic cycle of $X$ 
in $M$ equals $(-1)^{\dim X}[\Proj(\Sym_{\mathcal O_X}(\mathcal J_{Y,X}))]$.
\end{corol}

\begin{proof}
Immediate consequence of Theorem~\ref{qSym}: if $X$ satisfies the \x-condition, 
then by Proposition~\ref{symqsym} the algebras $\Sym_{\mathcal 
O_X}(\mathcal J_{Y,X})$ and $\qSym_{\mathcal O_X}(\mathcal J_{Y,X})$
are isomorphic in high degree, so they have the same $\Proj$.
\end{proof}

\subsection{The \x-files}\label{xfiles}
The `\x' in {\em \x-condition\/} has been chosen as it reminds us of the 
prototypical singularities satisfying it: the conic $xy=0$ in the plane is 
(i) a hypersurface with a nonsingular singularity subscheme; (ii) a 
hypersurface with quasi-homogeneous isolated singularities; and (iii) a 
divisor with normal crossing divisor. Each of these classes of 
hypersurfaces satisfies the \x-condition. In fact, as the interested 
reader may verify, in each of these cases the embedding of the singularity 
subscheme in the ambient space is `linear'.

Recall (\cite{MR93c:14007}) that an embedding of schemes $S\subset T$ is 
{\em linear\/} if the Rees algebra and the symmetric algebra of the ideal of 
$S$ in $T$ are isomorphic; it is {\em weakly linear\/} if the Rees algebra 
and the symmetric algebra are isomorphic in high degree, that is, if
$\Proj(\Sym_T(S))$ is isomorphic to the (Rees) blow-up of $T$ along
$S$.

\begin{prop}\label{weaklin}
Let $X$ be a hypersurface in a nonsingular variety $M$, with singularity 
subscheme~$Y$. If the embedding of $Y$ in $M$ is weakly linear, then
$X$ satisfies the \x-condition.
\end{prop}

\begin{proof}
With the notations in the proof of Proposition~\ref{symqsym}, the
embedding of $Y$ in $M$ is weakly linear if and only if $\Tors_d=0$
for $d\gg 0$. This implies $\Disc_d=0$ for $d\gg 0$, which is
equivalent to the \x-condition, as observed in that proof.
\end{proof}

For example, this implies immediately the \x-condition for the first case 
listed above: if $Y$ is nonsingular, then its embedding in $M$ is regular, 
hence linear, hence weakly-linear. However, we should remind the reader 
that the requirement that the singularity {\em subscheme\/} of a 
hypersurface be nonsingular is very strong; substantially stronger, for 
example, than the requirement that the singularity {\em locus\/} be 
nonsingular. Some constraints on this situation are studied in 
\cite{MR97b:14057}, \S3. Hypersurfaces whose singularity subscheme is 
nonsingular are in particular {\em nice\/} in the sense of 
\cite{math.AG/0107185}.

\begin{example} The plane curve $x^4+x^3 y^2+y^6=0$ has an isolated
singularity at the origin; the embedding of its singularity subscheme
in the plane is not linear.
\end{example}

This is checked by explicit calculations, which we performed with {\tt
Macaulay 2}.

It can also be shown that if the \x-condition implies that every vector 
tangent at a point $x$ to a stratum in a Whitney stratification of $X$ 
extends to fiberwise linear functions on $\mathcal P^1_M\mathcal O(X)$, 
tangent to nearby `level hypersurfaces'. In this sense, the \x-condition 
may be viewed as a strong regularity requirement on extensions of tangent 
vectors near strata of a Whitney stratification of $X$. It is known that 
tangent vectors to strata of a Whitney stratification of $X$ are suitably 
`close' to the tangent spaces of nearby level hypersurfaces (the so-called 
{\em $w_f$-condition of Thom}). This suggests that the techniques in 
\cite{MR94i:32059} or \cite{MR95h:32043} may be apt to characterizing 
hypersurfaces satisfying the \x-condition.

\subsection{}\label{nonred} 
For simplicity we have assumed that the hypersurface $X$ is reduced in the 
preceding subsections. It should be noted, however, that the 
quasi-symmetric blow-up is defined, and determines a cycle in $\Pbb(T^*M)$, 
regardless of whether $X$ is reduced or not. The arguments given above can be 
traced in this case, and show that this cycle is nothing but the characteristic 
cycle of the support $X_{\text{red}}$. This rather remarkable fact implies 
that simply setting $\csm(X):=\csm(X_{\text{red}})$ leads to a consistent 
theory of Chern-Schwartz-MacPherson classes, at least when $X$ is a 
hypersurface. 

We leave the details to the interested reader (cf.~\S2.1 in 
\cite{MR2001i:14009}).


\section{Shadows of blow-up algebras}\label{charcla}

\subsection{}\label{promise} 
Theorems~\ref{rees} and \ref{qSym} give intrinsic constructions of the two 
key cycles associated with~$X$. We would like to deal with the corresponding 
schemes $B\ell_YX$, $\qBl_YX$ as stand-alone entities, and determine 
precisely what type of information they carry in relation with the ambient 
nonsingular variety~$M$.

With this in mind, we first discuss the transformation $\mathcal L\leadsto 
\mathcal A$ mentioned in \S\ref{sabbah}, which produces the Chern-Mather, 
resp.~Chern-Schwartz-MacPherson classes from the conormal,
resp.~characteristic cycle; then we separate the r\^ole of the ambient
variety in this computation from that of the blow-ups themselves, and
find that the blow-ups carry `normal data' regarding the embedding
$X\subset M$.  This point of view unifies the computation of the
Chern-Mather and Chern-Schwartz-MacPherson classes with the approach
yielding the classes defined by Fulton and Fulton-Johnson
(\cite{MR85k:14004}, \cite{MR82c:14017}, and cf.~\S\ref{ffj} below).

\subsection{} 
If $\mathcal E$ is a locally free sheaf of rank $e+1$ on a scheme $S$, there 
is a precise structure theorem for the Chow group of the projective bundle 
$$\xymatrix{
{\Pbb(\mathcal E):=\Proj(\Sym \mathcal E^\vee)} \ar[r]^-\epsilon & S
}$$ 
(\cite{MR85k:14004}, \S3.3): every class $C\in A_r\Pbb(\mathcal E)$
can be written uniquely as
$$C=\sum_{j=0}^e c_1(\mathcal O(1))^j\cap \epsilon^*(C_{r-e+j})$$
where $\mathcal O(1)$ denotes the tautological line bundle on 
$\Pbb(\mathcal E)$, and $C_{r-e+j}\in A_{r-e+j}S$.

Therefore, knowledge of $C$ is equivalent to knowledge of the
collection of $e+1$ classes $C_{r-e},\dots,C_r$ on $S$.

\begin{defin}\label{shaddef}
We say that the class $C_{r-e}+\dots+C_r\in AS$ is the {\em shadow\/} of 
the class $C$.
\end{defin}

As its real world namesake, the shadow neglects some of the
information carried by the object that casts it. For example,
$c_1(\mathcal O(1))^j\cdot [\Pbb(\mathcal E)]$ has shadow $[S]$ for
all $j=0,\dots,e$.  However, a pure-dimensional class $C$ can be
reconstructed from its shadow if its dimension is known, as follows
immediately from the structure theorem recalled above.

It will be convenient to have a direct way to obtain the shadow of a
given class.

\begin{lemma}\label{shadowc}
The shadow of $C$ is the class
$$c(\mathcal E)\cap \epsilon_*\left(c(\mathcal O(-1))^{-1}\cap C\right)$$
\end{lemma}

\begin{proof}
Writing $C$ as above, we have
$$\aligned
c(\mathcal E)\cap \epsilon_*\left(c(\mathcal O(-1))^{-1}\cap C\right)
&=c(\mathcal E)\cap \epsilon_*\left(c(\mathcal O(-1))^{-1}\cap
\sum_{j=0}^e c_1(\mathcal O(1))^j\cap \epsilon^*(C_{r-e+j})\right)\\
&=\sum_{j=0}^e c(\mathcal E)\cap \epsilon_*\left( \sum_{k\ge j}
c_1(\mathcal O(1))^k\cap \epsilon^*(C_{r-e+j})\right)\quad.
\endaligned$$
Since $c_1(\mathcal O(1))^k\cap \epsilon^*\alpha=0$ for $0\le k<e$ and
any $\alpha\in A_*S$, this says
$$c(\mathcal E)\cap \epsilon_*\left(c(\mathcal O(-1))^{-1}\cap C\right)
=\sum_{j=0}^e c(\mathcal E)\cap \epsilon_*\left(c(\mathcal
O(-1))^{-1}\cap \epsilon^*(C_{r-e+j})\right)\quad.$$
Finally, this equals $\sum_{j=0}^e C_{r-e+j}$ by \cite{MR85k:14004},
Example~3.3.3.\end{proof}

\subsection{}\label{sabbah2} As recalled in \S\ref{sabbah}, MacPherson's 
natural transformation $c_*$ can be expressed by a two-step procedure:
($\mathcal C \leadsto \mathcal L$) taking the characteristic cycle
$\Ch(\varphi)$ of a constructibile function $\varphi$, and ($\mathcal
L\leadsto \mathcal A$) extracting a rational equivalence class from
the characteristic cycle. As the natural habitat of Lagrangian
cycles is the projectivized {\em co\/}tangent bundle $\Pbb(T^*M)$, we
find it convenient to arrange things so as to obtain a class $\check
c_*(\varphi)$ differing from $c_*(\varphi)$ by the sign of the
components of odd dimension:
$$\{\check c_*(\varphi)\}_r=(-1)^r \{c_*(\varphi)\}_r$$
in dimension $r$. For example,
$$\check c_*(\one_M)=(-1)^{\dim M} c(T^*M)\cap [M]$$ 
for the nonsingular ambient $M$.

\begin{lemma}\label{shad}
The class $\check c_*(\varphi)$ is the shadow of the characteristic cycle
$\Ch(\varphi)$.
\end{lemma}

\begin{proof} This is formula (12) on p.~67 of \cite{MR1795550},
filtered through Lemma~\ref{shadowc}. As observed in \cite{MR1795550},
this is in agreement with \cite{MR50:13587}.
\end{proof}

The statement of Lemma~\ref{shad}, while implicit in the existing
literature, is mysteriously absent in this explicit form relating the
transformation $\mathcal L\leadsto \mathcal A$ to the structure
theorem of the Chow group of projective bundles. This interpretation
streamlines the proof that $\mathcal L\leadsto \mathcal A$ is a
natural transformation; J\"org Sch\"urmann has independently made the
same observation \cite{schuer}.

\subsection{} We are ready to justify the title of this article. Denote 
by $\ccma(X)$, $\ccsm(X)$ respectively the classes obtained by
changing the sign of the components of odd dimension in $\cma(X)$,
$\csm(X)$.

\begin{theorem}\label{shadtheorem}
Let $X$ be a hypersurface of a nonsingular variety $M$, and let $Y$ be
its singularity subscheme. Then
\begin{itemize}
\item the shadow of $[B\ell_YX]$ is $(-1)^{\dim X}\ccma(X)$;
\item the shadow of $[\qBl_YX]$ is $(-1)^{\dim X}\ccsm(X)$.
\end{itemize}
\end{theorem}

\begin{proof}
This now follows from Theorems~\ref{rees} and \ref{qSym}, and 
Lemma~\ref{shad}.
\end{proof}

\subsection{} 
The next step in our program consists of carefully distinguishing the r\^ole 
of the ambient space and of the blow-ups in the statement of 
Theorem~\ref{shadtheorem}. There is an interesting twist to this story, 
which highlights the need for a subtle change of perspective.

We have so far focused on the ideal $\mathcal J_{Y,X}$ as the most natural 
source of information concerning the singularities of $X$; and indeed we 
have defined our main notions in \S\ref{qsymblowups} starting from the 
data of an ideal sheaf in $\mathcal O_X$. We are now going to shift the 
attention to a different coherent sheaf, defined for any subscheme $X$ of a 
nonsingular variety $M$; it will be easy to relate this sheaf to $\mathcal 
J_{Y,X}$ when $X$ is a hypersurface, and this will naturally extend 
quasi-symmetric blow-up algebras to this coherent sheaf. To summarize what 
we will find, this new algebras agree locally with the algebras obtained 
for $\mathcal J_{Y,X}$; in fact, their Proj will be isomorphic as schemes 
to the quasi-symmetric blow-ups of $\mathcal J_{Y,X}$. But the {\em 
algebras\/} carry more information than the schemes: the grading determines 
a line bundle on the blow-ups, and this information will turn out to be 
essential.

The new blow-up algebras will thus determine a Segre-class type of 
invariant, and we will show that using this invariant yields the 
Mather and Schwartz-MacPherson classes in essentially the same way as 
ordinary Segre classes of coherent sheaves, resp.~of cones lead to 
Fulton-Johnson, resp.~Fulton classes.

\subsection{}\label{ffj}
Here is a quick reminder concerning these latter two classes, in order to 
clarify the context underlining our motivation.

If $Z$ is any scheme embedded in a nonsingular variety $M$ (of 
dimension $>\dim Z$ for convenience), there are several ways to obtain `normal 
data' relating to the embedding. For example, such data is carried by the 
conormal sheaf $\mathcal N_ZM=\mathcal J_{Z,M}\otimes_{\mathcal O_M} \mathcal 
O_Z$, and can be effectively encoded in the {\em Segre class\/} of this 
coherent sheaf, defined by
$$s(\mathcal N_ZM):=p_*\sum c(\mathcal O(1))^i\cap [\Proj(\Sym_{\mathcal O_Z}
(\mathcal N_ZM))]$$
where $p$ is the structure morphism on $\Proj$.

The class $c(TM)\cap s(\mathcal N_ZM)$ is the {\em Fulton-Johnson\/} class 
of $Z$. It can be shown to be independent of the embedding, and agrees with 
the total Chern class of the tangent bundle of $Z$ when $Z$ is nonsingular
(cf.~\cite{MR82c:14017} or \cite{MR85k:14004}, Example~4.2.6 (c)).

A different way to access normal data amounts to taking a Rees point of 
view rather than a Sym point of view. Replacing
$$\Sym_{\mathcal O_Z}(\mathcal N_ZM)=\Sym_{\mathcal O_M}(\mathcal J_{Z,M}) 
\otimes_{\mathcal O_M} \mathcal O_Z$$
by 
$$\Rees_{\mathcal O_M}(\mathcal J_{Z,M}) \otimes_{\mathcal O_M} 
\mathcal O_Z$$
defines the {\em normal cone\/} of $Z$ in $M$, whose Segre class (again 
defined by pushing forward powers of the first Chern class of $\mathcal 
O(1)$) is properly called the {\em Segre class of $Z$ in $M$,\/} $s(Z,M)$.

Applying the same principle as above leads to defining the class 
$c(TM)\cap s(Z,M)$, which can again be shown to be independent of the 
embedding (cf.~\cite{MR85k:14004}, Example~4.2.6), and which again agrees 
with the total Chern class of the tangent bundle of $Z$ when $Z$ is 
nonsingular. This class is called the {\em Fulton\/} class of $Z$.

\subsection{}
How else can one extract normal data from an embedding $Z\subset M$ of a 
scheme in a nonsingular variety? Again we assume that $\dim M>\dim Z$. There 
is a surjection
$$\Omega^1_M|_Z \longrightarrow \Omega^1_Z \longrightarrow 0\quad,$$
from which we obtain the exact sequence
$$0 \longrightarrow \cHom(\Omega^1_Z,\mathcal O_Z) \longrightarrow 
\cHom(\Omega^1_M|_Z,\mathcal O_Z) \longrightarrow \mathcal T_ZM 
\longrightarrow 0\quad,$$
defining the coherent sheaf $\mathcal T_ZM$ on $Z$. If $Z$ is
nonsingular then $\mathcal T_ZM$ is locally free, and in fact it is
the sheaf of sections of the normal bundle of $Z$ in $M$.

Now our idea consists of following the same guiding principle which rules 
in \S\ref{ffj}, but employing Segre classes obtained from quasi-symmetric 
algebras associated with $\mathcal T_ZM$. As things stand now, we only 
have defined such objects for ideals, and this limits the scope of our 
aim. However, in the case we have considered in \S\ref{cycles} and in 
Theorem~\ref{shadtheorem} the day is saved by a special form taken by 
$\mathcal T_ZM$.

\begin{lemma}\label{txm}
If $Z=X$ is a hypersurface in a nonsingular variety $M$, with line 
bundle $\mathcal L$ and singularity subscheme $Y$, then 
$\mathcal T_XM=\mathcal J_{Y,X}\otimes_{\mathcal O_X}\mathcal L$.
\end{lemma}

\begin{proof}
Let $\mathcal J=\mathcal J_{X,M}$ denote the ideal of $X$ in $M$. Taking
$\cHom(-,\mathcal O_X)$ in the exact sequence of differential gives the 
exact sequence
$$0 \longrightarrow \cHom(\Omega^1_X,\mathcal O_X) \longrightarrow 
\cHom(\Omega^1_M|_X,\mathcal O_X) \longrightarrow \cHom(\mathcal J/\mathcal 
J^2,\mathcal O_X)\quad.$$
A local computation determines the image of the rightmost map as 
the subsheaf of $\mathcal O_X$-morphisms $\mathcal J/\mathcal J^2
 \longrightarrow \mathcal O_X$ factoring through $\mathcal J_{Y,X}$. In other 
words, if $X$ is a hypersurface in $M$ then
$$\mathcal T_XM=\cHom(\mathcal J/\mathcal J^2,\mathcal J_{Y,X}) = 
\mathcal L|_X \otimes_{\mathcal O_X} \mathcal J_{Y,X}\quad,$$
as claimed.
\end{proof}

By virtue of Lemma~\ref{txm} we can make sense of quasi-symmetric algebras 
of $\mathcal T_XM$ if $X$ is a hypersurface in $M$. The two extremes in the 
range of quasi-symmetric algebras are the following two definitions:
$$\aligned
\qSym_{X\subset X}(\mathcal T_XM) &:= \Rees_{\mathcal O_X}(\mathcal 
J_{Y,X})\otimes_{\mathcal O_X}\mathcal L\\
\qSym_{X\subset M}(\mathcal T_XM) &:= \qSym_{\mathcal O_X}(\mathcal 
J_{Y,X})\otimes_{\mathcal O_X}\mathcal L
\endaligned$$
and the corresponding Segre class-like notions:
$$\aligned
\csma(X,M) &:=p_* \sum c(\mathcal O(1))^i \cap [\Proj(\qSym_{X\subset X}
(\mathcal T_XM))] \\
\csma(X,M) &:=p_* \sum c(\mathcal O(1))^i \cap [\Proj(\qSym_{X\subset M}
(\mathcal T_XM))]
\endaligned$$
where $p$ denotes the projection from the corresponding Proj, and 
$\mathcal O(1)$ is the tautological line bundle. We remark that the two 
Proj equal $B\ell_Y M$, $\qBl_Y M$ as schemes---only the tautological 
bundles are affected upon tensoring by $\mathcal L$.

We have defined a `checked' notion of Segre class in view of artificially 
taking a dual that brings us back to the {\em tangent\/} world. So we set
$$\aligned
\sma(X,M) &:= (-1)^{\dim X}\sum_{r\ge 0} (-1)^r \csma(X,M)_r \\
\ssm(X,M) &:= (-1)^{\dim X}\sum_{r\ge 0} (-1)^r \cssm(X,M)_r
\endaligned$$
where subscripts mark dimensions; that is, we change the sign of components 
in the checked Segre classes of every other {\em co\/}dimension in $X$.

\begin{example}
If $X$ is a nonsingular hypersurface, then all notions of Segre class
coincide: $s(\mathcal N_XM)=s(X,M)=\sma(X,M)=\ssm(X,M)=c(\mathcal
L)^{-1}\cap [X]$.  If $X$ may be singular, but satisfies the
\x-condition (see \S\ref{xcond}), then $\ssm(X,M)=s(\mathcal
J_{Y,X}\otimes\mathcal L)$.
\end{example}

\subsection{}
Summarizing, we have extracted normal data from our hypersurface $X$ in $M$ 
by defining a coherent sheaf $\mathcal T_XM$ in a rather simple-minded way 
from the exact sequence of differentials of $X$; adapting to $\mathcal 
T_XM$ the construction of \S\ref{qsymblowups}; and defining from the 
resulting blow-up algebra a notion of Segre class. These classes achieve 
precisely what we set out to do, that is, they yield the Chern-Mather and 
Chern-Schwartz-MacPherson classes by the same method behind the classes of 
Fulton and Fulton-Johnson. That is:

\begin{theorem}\label{lean}
Let $X$ be a hypersurface in a nonsingular variety $M$. Then
\begin{itemize}
\item $\cma(X)=c(TM)\cap \sma(X,M)$;
\item $\csm(X)=c(TM)\cap \ssm(X,M)$.
\end{itemize}
\end{theorem}

\begin{proof}
We will give the argument for the second equality; the first is treated 
similarly.

Tensoring by $\mathcal L$ the epimorphism
$$\xymatrix{
\Sym((\Omega^1_M\otimes \mathcal L)^\vee|_X) \ar@{->>}[r] &
\qSym_{\mathcal O_X}(\mathcal J_{Y,X})
}$$
from \S\ref{properplace} we obtain 
$$\xymatrix{
\Sym((\Omega^1_M)^\vee|_X) \ar@{->>}[r] & 
\qSym_{X\subset M}(\mathcal T_XM)\quad,
}$$
inducing the embedding
$$\qBl_YX \hookrightarrow \Pbb(T^*M)$$
realizing the characteristic cycle of $X$ (by Theorem~\ref{qSym}), and 
showing that the restriction of $\mathcal O(-1)$ to $\qBl_YX$ is the universal 
bundle $\mathcal O(-1)$ of $\Proj(\qSym_{X\subset M}(\mathcal T_XM))$. By 
Lemma~\ref{shadowc}, the shadow of the blow-up algebra $\qBl_YX$ is computed by
$$c(T^*M)\cap (c(\mathcal O(-1))^{-1}\cap [\qBl_YX]) = c(T^*M)\cap 
\cssm(X,M)\quad.$$
This equals $(-1)^{\dim X}\ccsm(X)$, by Theorem~\ref{shadtheorem}. The 
equality for $\csm(X)$ follows by changing the sign of the components of 
every other codimension.\end{proof}

\subsection{} At this point it is only too natural to pose the problem 
of defining quasi-symmetric algebras for coherent sheaves so as to 
validate Theorem~\ref{lean} for more general schemes $X$, following the same 
strategy (that is, by obtaining Segre classes from the quasi-symmetric 
algebras of $\mathcal T_XM$). The advantage in formulas such as those in 
Theorem~\ref{lean} is not only theoretical: these formulas can be 
implemented in procedures for symbolic computation programs such as {\tt 
Macaulay2.\/} At present a routine is implemented that computes 
Chern-Schwartz-MacPherson classes of projective schemes
(\cite{math.AG/0204230}), exploiting the hypersurface case in order
to compute classes in the general case, by a computationally expensive
`inclusion-exclusion' procedure.

An upgrade of Theorem~\ref{lean} to more general schemes would bring 
about a drastic improvement in the speed of such routines.

Regarding a possible definition of quasi-symmetric algebras for coherent 
sheaves, this would presumably pivot on a good notion of Rees algebra of 
a module; such notions have been introduced and studied, primarily by 
Artibano Micali (starting with \cite{MR31:1275}). Even in the simpler case 
of ideals treated here, it would be quite interesting to relate our 
construction with the ideals defined by Micali in loc.cit., 
interpolating between the symmetric and the Rees algebras.



\end{document}